\documentclass[11pt,reqno]{amsart}

\usepackage{newtxtext}
\usepackage{graphicx}
\usepackage{newtxmath}

\usepackage{cite}

\usepackage{url}
\usepackage{enumitem}

\newtheorem{theorem}{Theorem}

\newtheorem{lemma}{Lemma}
\newtheorem{corollary}{Corollary}
\newtheorem{notation}{Notation}
\newtheorem{definition}{Definition}

\newtheorem{remark}{Remark}

\usepackage[pdftex,unicode,bookmarksnumbered]{hyperref} 

\hypersetup{
 pdftitle={On $\theta$-extension of continuous mapping},
 pdfauthor={Andrew Ryabikov},
 pdfsubject={Extension of maps},
 pdfkeywords={
  mapping extension problem{,}
  topological space{,}
  weakly continuous mapping{,}
  $\theta$-continuous mapping
 },
}

\hypersetup{
 colorlinks=true,
 allcolors=[rgb]{0,0.5,0.5},
 urlcolor=[rgb]{0.58,0.0,0.83},
 breaklinks=true,
}

\begin{document}     


\title{On $\theta$-extension of continuous mapping}

\author{Andrew Ryabikov}

\address{
Federal Research Center "Computer Science and Control" of the Russian Academy of Sciences \\
Vavilova str. 44\\ 
119933 Moscow\\ 
Russia
}

\email{ariabikov@gmail.com}



\begin{abstract}
We consider the problem of constructing a weakly-continuous mapping extending continuous mapping defined on a dense set of a topological space to the entire space. Theorem on necessary and sufficient conditions for the existence of such an extension are proved. The axioms of separation and compactness of spaces are not assumed.
\end{abstract}

\maketitle


\section*{Introduction}

Let a mapping $f$ of an everywhere dense set $S$ of some topological space $X$ into a topological space $Y$ be given. The problem of extending the mapping $f$ is understood as finding a continuous (weakly continuous) mapping $F:X\rightarrow Y$ whose restriction to the set $S$ coincides with the mapping $f$.

\par
The image $F$ of some point $x$ not lying in $S$ can be obtained using a centered system of neighbourhoods in which the original mapping $f$ is defined. For example, in \cite{Bourb39,Vul52,Taim52,Vel65,Engel64,Osipov17,Osm86} $F(x)$ is defined as a common element of closed sets $[f(S\cap P)]$, where $P\in \mathcal{N}(x)$, and $\mathcal{N}(x)$ is a local base of the topological space $X$. If the separation axioms are satisfied and some restrictions are imposed on the original mapping, we obtain a continuous mapping $F$.
\par
The question arises of how to define the desired mapping if the spaces $X$ and $Y$ are not compact and the fulfilment of the separation axioms is not assumed. We know that a topological space is defined by a system of open sets. Having defined such a system, we can speak, for example, about the existence or absence of limit points of sets. On the other hand, the limit points of the sets themselves can characterize the topological space. This paper is based on the idea that a continuous mapping from $X$ to $Y$ connects the limit points of sets in $X$ with the limit points of their images in $Y$. In this paper, we prove a theorem on necessary and sufficient conditions for the existence of $\theta$-continuous mapping that extends the original mapping $f$ to the entire space $X$.


\section*{Prerequisites}

\begin{notation}
We will denote the closure of a set by square brackets $[\,]$. By a neighbourhood of a point in a topological space we will mean an open set containing the given point. Let $\alpha>0$ be an ordinal.
\par
The definition of $\theta_\alpha$ - continuous mapping and $\alpha$-hull of a set is given in \cite{Osipov17}. We present them below.
\end{notation}

\begin{definition}
A neighbourhood $U_{\alpha}$ of a set $A\subseteq Y$ is called an $\alpha$-hull of $A$ if there exists a collection $\{U_{\beta}\}_{\beta\leq \alpha}$ of open sets containing $A$ for which $[U_{\beta}]\subseteq U_{\beta+1}$ holds for all $\beta+1\leq\alpha$ and $U_{\alpha}=\bigcup\limits_{\beta\leq \alpha} U_{\beta}$.
\end{definition}

\begin{definition}
The set $[A]_{\theta_\alpha}$ is called the $\theta_\alpha$-closure of the set $A$. Point $x\in [A]_{\theta_\alpha}$ if the $[U_{\alpha}(x)]\bigcap A \neq \emptyset$ for any $\alpha$-neighbourhood $U_{\alpha}(x)$.
\end{definition}

\begin{definition}
A mapping $F:X\rightarrow Y$ is called $\theta_\alpha$ - continuous if for any point $x\in X$ and for any $\alpha$-hull $U_{\alpha}$ of the image $F(x)$ there is a neighbourhood $Ox$ such that $F(Ox)\subseteq [U_{\alpha}]$.
\end{definition}
In the case $\alpha=1$ we get the definition of a weakly continuous mapping or $\theta$-continuous mapping. 
\par 
For $\alpha=0$ we set $[A]_{\theta_\alpha} = [A]$, and $\theta_\alpha$-continuous mapping is considered a continuous mapping.


\section*{Mapping Extension Theorem}

Let $X$ be a topological space and $S\subset X$ be a dense set. Let $f:S\rightarrow Y$ be a continuous mapping. Denote by $K(x)$ the set of all subsets of $S$ whose closure contains $x$, i.e. $K(x)=\{ M\subseteq S | x\in [M]\}$.
\par Let us formulate a theorem on the $\theta$-continuous extension of the mapping $f$ to the entire space $X$.
\begin{theorem}
For the existence of a $\theta$-continuous mapping $F:X\rightarrow Y$ whose restriction to $S$ coincides with $f$, it is necessary that $\bigcap\limits_{M\in K(x)}[f(M)]_\theta \neq \emptyset$ for any $x\in X$ and it is sufficient that $\bigcap\limits_{M\in K(x)}[f(M)]\neq \emptyset$ for any $x\in X$.
\end{theorem}

To prove Theorem 1, we need an auxiliary lemma, which we formulate for $\theta_\alpha$ - continuous mappings.

\begin{lemma}
A mapping $F:X\rightarrow Y$ is $\theta_{\alpha}$-continuous if and only if $F([A]) \subseteq [F(A)]_{\theta_\alpha}$ for any $A\subset X$.
\end{lemma}

\begin{proof}
{\it Necessity.} Let $F:X\rightarrow Y$ be $\theta_{\alpha}$-continuous. We will show that $F([A]) \subseteq [F(A)]_{\theta_\alpha}$ for any $A\subset X$.
\par
Let this not be so and there is a set $A\subset X$ and a point $x\in [A]$ such that $F(x)\notin [F(A)]_{\theta_\alpha}$. Let $\alpha>0$. Then there exists a neighbourhood $U_{\alpha}$ of $F(x)$ such that $[U_{\alpha}]\bigcap F(A) = \emptyset$. Since $F$ is $\theta_{\alpha}$-continuous, we have $F(Ox)\subseteq [U_{\alpha}]$. Consequently, $F(Ox)\bigcap F(A) = \emptyset$. But $Ox\bigcap A \neq \emptyset$, and hence $F(Ox)\bigcap F(A) \neq \emptyset$. We have arrived at a contradiction.
\par
{\it Sufficiency.} Let $F([A]) \in [F(A)]_{\theta_\alpha}$ for any $A\subset X$. We will show that the mapping $F:X\rightarrow Y$ is $\theta_{\alpha}$-continuous.
\par
Suppose that $F$ is not $\theta_{\alpha}$-continuous. Let $\alpha>0$. Then for some neighbourhood $U_{\alpha}$ of $F(x)$ there exists $A\subset X: x\in [A]$ such that $[U_{\alpha}]\bigcap F(A) = \emptyset$. By the hypothesis of the theorem, $F(x) \in F([A]) \subseteq [F(A)]_{\theta_\alpha}$. Therefore, $[U_{\alpha}]\bigcap F(A) \neq \emptyset$. We have arrived at a contradiction. Thus, $F$ is $\theta_{\alpha}$-continuous. If $\alpha=0$ it is reasonable to replace $[U_{\alpha}]$ with some neighbourhood $U$ of the point $F(x)$. The lemma is proved.

\end{proof}

Let us proceed to the proof of Theorem 1.

\begin{proof}
{\it Necessity.}
Let there exist a $\theta$-continuous mapping $F:X\rightarrow Y$ whose restriction to $S$ coincides with $f$. By Lemma 1, $F(x) \in F([M])\subseteq [F(M)]_\theta$ for all $M\in K(x)$. Whence the assertion of the theorem follows.

{\it Sufficiency.}
We define the mapping $F$ as follows. Let $F(x)=f(x)$ for all $x\in S$. For $x\in X\setminus S$, we take as $F(x)$ an arbitrary common point of the sets $[f(M)]$, where $M\in K(x)$.
As it follows from the Lemma 1 for $x\in S$ the mapping $F(x)$ is a common point of $[f(M)]$, where $M\in K(x)$. Let us prove that the mapping $F$ is $\theta$-continuous.
\par
Take an arbitrary neighbourhood $OF(x)$ of the image $F(x)$. Consider the operator $C(M) = \{z\in M| f(z)\in OF(x)\}$, where $M\in K(x)$, which restricts the set $M\in K(x)$ according to a given rule. Let us show that $C(M)\in K(x)$.
\par
Assume the opposite. Then there exists a neighbourhood $Ox$ such that \par\noindent $Ox\bigcap C(M) = \emptyset$. But, on the other hand, $Ox\bigcap M \in K(x)$, and therefore, by the hypothesis of the theorem, $F(x)\in [f(Ox\bigcap M)]$. Consequently, there exists a point $z\in Ox\bigcap M$ at which $f(z)\in OF(x)$. We have arrived at a contradiction with the fact that $Ox\bigcap C(M) = \emptyset$. Thus, $C(M)\in K(x)$ for any $M\in K(x)$.
\par
We will show that there is a neighbourhood $Ox : Ox\subseteq [C(S)]$. Suppose that this is not so. Then there exists a set $A\subset X : x\in [A]$ and $A\bigcap [C(S)]=\emptyset$. We denote $OA = X\setminus [C(S)]$. Since $x\in [OA]$, then $M = OA\bigcap S \in K(x)$. Clearly, $M\bigcap [C(S)]=\emptyset$.
\par
The following embeddings are true: $C(M)\subseteq M$ and $C(M)\subseteq C(S)$. Thus, $M\bigcap C(S) \neq \emptyset$, which contradicts $M\bigcap [C(S)]=\emptyset$, and hence there exists a neighbourhood $Ox : Ox\subseteq [C(S)]$.
\par
Consider $x_0 \in Ox \subseteq [C(S)]$. Clearly, $[C(S)]\in K(x_0)$. By the hypothesis of the theorem, $F(x_0)\in [f(C(S))]$. Therefore, $F(x_0)\in [OF(x)]$. Thus, $F(Ox)\subseteq [OF(x)]$. Theorem 1 is proved.
\end{proof}

\begin{remark}
The requirement $F|_S = f$ is not necessary. We only need the condition $F(x)\in [f(M)]$ for all $M\in K(x)$. In this case, the $\theta$-continuous mapping $F$ approximates $f$. 
\end{remark}

\begin{corollary}
For the existence of a continuous mapping $F$ of a space $X$ into a regular space $Y$ whose restriction to $S$ coincides with $f$, it is necessary and sufficient that $\bigcap\limits_{M\in K(x)}[f(M)]\neq \emptyset$ for any $x\in X$.
\end{corollary}

\newpage

Let us give a simple example showing that the sufficient condition of Theorem 1 is not necessary.

\par\medskip\noindent
{\bf Example.} As $X$ we take the segment $[0,1]$ with the Euclidean metric. Let $S=X\setminus\{0.5\}$. On the set $Y=\{0,1,2\}$ we define the topology as shown in Fig. 1, where the neighborhoods are marked in blue. We define a continuous mapping $f:S\rightarrow Y$ as follows:
\begin{equation}\label{eq1}
f(x)=\left\{
\begin{array}{ll}
0, & \text{for } x\in [0, 0.5)\\
1, & \text{for } x\in (0.5, 1].
\end{array}
\right.
\end{equation}
\par
Let $F(x=0.5) = 1$, and let the restriction of $F$ to $S$ coincide with $f$. On the one hand, the mapping $f$ is continuous, and $F$ is a $\theta$-continuous extension of $f$ to $X$. On the other hand, the sufficient condition of Theorem 1 is not satisfied, since $\bigcap\limits_{M\in K(x)}[f(M)]=\emptyset$ for $x=0.5$.

\begin{figure}[ht] \centering
\includegraphics[height=0.2\textheight]{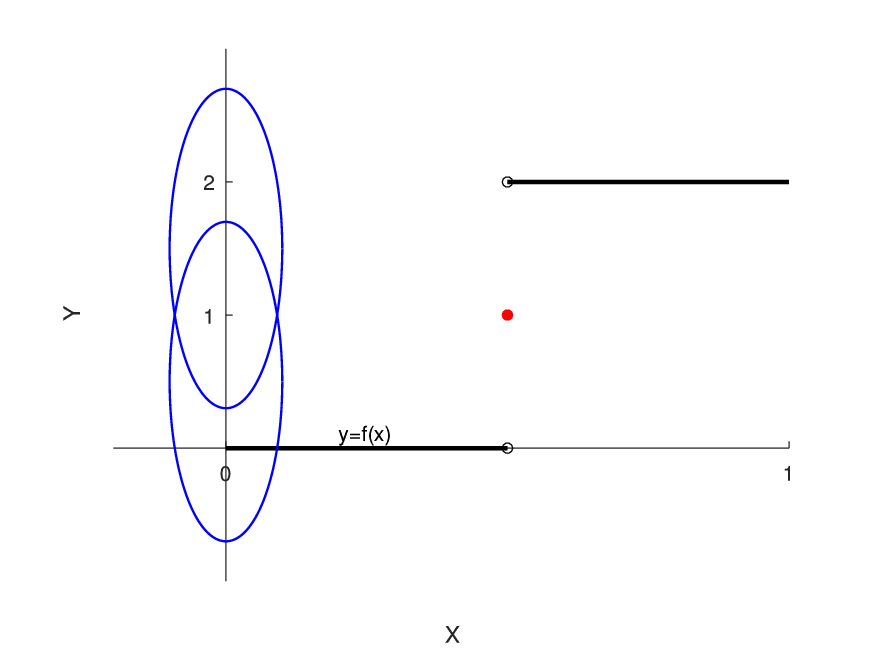}
\caption{Example of mappings $f$ and $F$}
\label{fig1}
\end{figure}

\section*{Conclusion}

The proved theorem justifies the existence of such a $\theta$-continuous mapping $F:X\rightarrow Y$ whose restriction to the set $S$ coincides with the mapping $f$. When constructing the mapping $F$, we assigned to each point $x$ of the space $X$ an arbitrary common point of contact of the sets $f(M)$, where $M\in K(x)$. It is possible that the choice of the desired $F(x)$ among the points of contact (the simplest limit points) of the sets $f(M)$ does not allow us to construct a continuous mapping $F$ that extends the original $f$. Further research will be aimed at choosing such a class of limit points within which we will be able to construct a continuous mapping $F$.




\begin{thebibliography}{1}
\bibitem{Bourb39}
N.~Bourbaki, J. Dieudonne, \textit{Notes de Teratopologie II}, Revue scientifique, 1939,  180--181 

\bibitem{Vul52}
B.\,Z.~Vulikh, \textit{On the extension of continuous functions in topological spaces},
Mat. Sb. \textbf{30} (1952), no. 1, 167--170  

\bibitem{Taim52} 
A.\,D.~Taimanov, \textit{On extension of continuous mappings of topological spaces},
Mat. Sb. \textbf{31} (1952), no. 2, 459--463


\bibitem{Vel65} 
N.\,V.~Velichko, \textit{O prodolgenii otobragenii topologicheskih prostranstv},
Sib. Mat. Zh. \textbf{6} (1965), no. 1, 64--69


\bibitem{Engel64} 
R.~Engelking, \textit{Remarks on real-compact spaces},
Fundam. Math. \textbf{55} (1964), 284--304

\bibitem{Osipov17} 
A.\,V.~Osipov, \textit{On extension functions for image space with different separation axioms},
Topology and its Applications \textbf{222} (2017), 70--76

\bibitem{Osm86} 
P.\,K.~Osmatesku, \textit{The extension of continuous maps},
Russian Mathematical Surveys \textbf{41} (1986), no. 6, 215--216


\bibitem{Oldair23} 
C.~Oldair-Renteria, \textit{On extension of fibrewise continuous and fibrewise almost-continuous functions},
Topology and its Applications \textbf{336} (2023)

\end{thebibliography}
\end{document}